# Null Similar Curves with Variable Transformations in Minkowski 3-space $E_1^3$

**Mehmet Önder**
*Celal Bayar University, Faculty of Science and Arts, Department of Mathematics, Muradiye Campus, 45047 Muradiye, Manisa, Turkey. E-mail: mehmet.onder@bayar.edu.tr*

**Abstract**
In this study, we define a family of null curves in Minkowski 3-space $E_1^3$ and called null similar curves. We obtain some properties of these special curves. We show that two null curves are null similar curves if and only if these curves form a null Bertrand pair. Moreover, we obtain that the family of null geodesics and null helices form the families of null similar curves with variable transformation.

## 1. Introduction

In the study of relativity theory, a unit speed future-pointing timelike curve in a space-time (a connected and time-oriented 4-dimensional Lorentz manifold) is thought as a locus of a material particle in the space-time. The unit-speed parameter of this curve is called the proper time of a material particle. Analogue to timelike curves, in relativity theory a future-pointing null geodesic is thought as the locus of a lightlike particle. More generally, from the differential geometric point of view, the study of null curves has its own geometric interest. Because the other curves (spacelike and timelike curves) of Lorentz space can be studied by a similar approach to that studied in positive definite Riemannian geometry. Moreover, null curves have different properties from spacelike and timelike cuvres and the results of null curve theory are not analogues to Riemannian case. So, motivated by the growing importance of null curves in mathematical physics, some special subjects of the curve theory have been studied for null curves by many mathematicians. Ferrandez, Gimenez and Lucas have defined and studied null helices in Lorentzian space forms [4]. Honda and Inoguchi have given a characterization of null cubics [5]. Later, Balgetir, Bektaş and Inoguchi have considered the notion of Bertrand curves for Cartan framed null curves and showed that null Bertrand curves are null geodesics or Cartan framed null curves with constant second curvature [1]. Moreover, Öztekin and Ergüt have given the characterizations of null Mannheim curves in $E_1^3$ [7].

Recently, a new definition of the special curves has been given by El-Sabbagh and Ali [3]. They have called these new curves as similar curves with variable transformation and defined as follows: Let $\psi_\alpha(s_\alpha)$ and $\psi_\beta(s_\beta)$ be two regular curves in $E^3$ parameterized by arc lengths $s_\alpha$ and $s_\beta$ with curvatures $\kappa_\alpha$, $\kappa_\beta$, torsions $\tau_\alpha$, $\tau_\beta$ and Frenet frames $\{\vec{T}_\alpha, \vec{N}_\alpha, \vec{B}_\alpha\}$ and $\{\vec{T}_\beta, \vec{N}_\beta, \vec{B}_\beta\}$, respectively. Then, $\psi_\alpha(s_\alpha)$ and $\psi_\beta(s_\beta)$ are called similar curves with variable transformation $\lambda_\beta^\alpha$ if there exists a variable transformation

$$s_\alpha = \int \lambda_\beta^\alpha(s_\beta) ds_\beta,$$

of the arc lengths such that the tangent vectors are the same for two curves, i.e., $\vec{T}_\alpha = \vec{T}_\beta$ for all corresponding values of parameters under the transformation $\lambda_\beta^\alpha$. All curves satisfying this condition is called a family of similar curves. Moreover, they have obtained some properties of the family of similar curves [3].

In this work we think the notion of similar curves for Cartan framed null curves in $E_1^3$. We give some theorems characterizing these special null curves and we show that the family of null geodesics and null helices form the families of null similar curves with variable transformation.

## 2. Preliminaries

Let $E_1^3$ be a Minkowski 3-space with natural Lorentz Metric

$$\langle,\rangle = -dx_1^2 + dx_2^2 + dx_3^2, \tag{1}$$

where $(x_1, x_2, x_3)$ is a rectangular coordinate system of $E_1^3$. According to this metric, in $E_1^3$ an arbitrary vector $\vec{v} = (v_1, v_2, v_3)$ can have one of three Lorentzian causal characters; it can be spacelike if $\langle \vec{v}, \vec{v} \rangle > 0$ or $\vec{v} = 0$, timelike if $\langle \vec{v}, \vec{v} \rangle < 0$ and null (lightlike) if $\langle \vec{v}, \vec{v} \rangle = 0$ and $\vec{v} \neq 0$ [6]. Similarly, an arbitrary curve $\vec{a} = \vec{a}(s)$ can locally be spacelike, timelike or null (lightlike), if all of its velocity vectors $\vec{a}'(s)$ are spacelike, timelike or null (lightlike), respectively. For the vectors $\vec{x} = (x_1, x_2, x_3)$ and $\vec{y} = (y_1, y_2, y_3)$ in $E_1^3$, the vector product of $\vec{x}$ and $\vec{y}$ is defined by

$$\vec{x} \wedge \vec{y} = \begin{vmatrix} e_1 & -e_2 & -e_3 \\ x_1 & x_2 & x_3 \\ y_1 & y_2 & y_3 \end{vmatrix} = (x_2 y_3 - x_3 y_2, x_1 y_3 - x_3 y_1, x_2 y_1 - x_1 y_2) \tag{2}$$

where

$$\delta_{ij} = \begin{cases} 1 & i = j, \\ 0 & i \neq j, \end{cases} \text{ and } e_i = (\delta_{i1}, \delta_{i2}, \delta_{i3}), \ e_1 \wedge e_2 = -e_3, \ e_2 \wedge e_3 = e_1, \ e_3 \wedge e_1 = -e_2. \tag{3}$$

**Definition 2.1.** A curve $\vec{a} = \vec{a}(s)$ in $E_1^3$ is called a null Frenet curve (or a Cartan-framed null curve) if it admits a Frenet frame field $\{\vec{\alpha}, \vec{\gamma}, \vec{\beta}\}$ such that

$$\begin{cases} \vec{\alpha}' = \kappa \vec{\beta}, \\ \vec{\gamma}' = \tau \vec{\beta}, \\ \vec{\beta}' = -\tau \vec{\alpha} - \kappa \vec{\gamma}, \end{cases} \tag{4}$$

with

$$\vec{\alpha} = \frac{d\vec{a}}{ds}, \ \langle \vec{\alpha}, \vec{\alpha} \rangle = \langle \vec{\gamma}, \vec{\gamma} \rangle = 0, \ \langle \vec{\beta}, \vec{\beta} \rangle = \langle \vec{\alpha}, \vec{\gamma} \rangle = 1 \tag{5}$$

where $\vec{\beta}$ is defined by $\vec{\beta} = \vec{\alpha} \times \vec{\gamma}$. The functions $\kappa$ and $\tau$ are called curvature and torsion of $\vec{a}(s)$, respectively. We called the vector fields $\vec{\alpha}, \vec{\gamma}, \vec{\beta}$ as tangent vector field, principal normal vector field and binormal vector field, respectively [5].

**Definition 2.2.** Let $\vec{a}(s_a)$ and $\vec{b}(s_b)$ be Cartan framed null curves in $E_1^3$ with frames $\{\vec{\alpha}_a, \vec{\gamma}_a, \vec{\beta}_a\}$ and $\{\vec{\alpha}_b, \vec{\gamma}_b, \vec{\beta}_b\}$, respectively. Then a pair of curves $(\vec{a}, \vec{b})$ is called a (null) Bertrand pair if the vector fields $\vec{\beta}_a$ and $\vec{\beta}_b$ are linearly dependent [1].

**Definition 2.3.** A null curve $\vec{a}(s_a)$ is said to be a null helix if it has constant Cartan curvatures [4].

After these definitions we give and prove the following theorem which will be used in the next section.

**Theorem 2.1.** *Let $\vec{a} = \vec{a}(s)$ be a Cartan framed null curve with arclength parameter $s$. Suppose that $\vec{a} = \vec{a}(\varphi)$ is an another parametrization of the curve by the parameter $\varphi(s) = \int \kappa(s) ds$. Then the null tangent vector field $\vec{\alpha}$ satisfies a vector differential equation of third order given by*

$$\frac{d^3 \vec{\alpha}}{d\varphi^3} + 2f(\varphi) \frac{d\vec{\alpha}}{d\varphi} + \frac{df(\varphi)}{d\varphi} \vec{\alpha} = 0, \tag{6}$$

*where $f(\varphi) = \dfrac{\tau(\varphi)}{\kappa(\varphi)}$.*

**Proof:** If we write derivatives given in (4) according to $\varphi$, we have

$$\begin{cases} \dfrac{d\vec{\alpha}}{d\varphi} = \dfrac{d\vec{\alpha}}{ds} \dfrac{ds}{d\varphi} = (\kappa \vec{\beta}) \dfrac{1}{\kappa} = \vec{\beta}, \\[6pt] \dfrac{d\vec{\gamma}}{d\varphi} = \dfrac{d\vec{\gamma}}{ds} \dfrac{ds}{d\varphi} = (\tau \vec{\beta}) \dfrac{1}{\kappa} = f(\varphi) \vec{\beta}, \\[6pt] \dfrac{d\vec{\beta}}{d\varphi} = \dfrac{d\vec{\beta}}{ds} \dfrac{ds}{d\varphi} = (-\tau \vec{\alpha} - \kappa \vec{\gamma}) \dfrac{1}{\kappa} = -f(\varphi) \vec{\alpha} - \vec{\gamma} \end{cases} \tag{7}$$

respectively, where $f(\varphi) = \dfrac{\tau(\varphi)}{\kappa(\varphi)}$. Then corresponding matrix form of (4) can be given by

$$\begin{bmatrix} d\vec{\alpha}/d\varphi \\ d\vec{\gamma}/d\varphi \\ d\vec{\beta}/d\varphi \end{bmatrix} = \begin{bmatrix} 0 & 1 & 0 \\ 0 & f(\varphi) & 0 \\ -f(\varphi) & 0 & -1 \end{bmatrix} \begin{bmatrix} \vec{\alpha} \\ \vec{\beta} \\ \vec{\gamma} \end{bmatrix}. \tag{8}$$

From the first and third equations of new Frenet derivatives (8) we have

$$\vec{\gamma} = -\left( \frac{d^2 \vec{\alpha}}{d\varphi^2} + f(\varphi) \vec{\alpha} \right). \tag{9}$$

Substituting the above equation in the second equation of (8) we have desired equation (6).

### 3. Null Similar Curves with Variable Transformation in $E_1^3$

In this section we introduce the definition and characterizations of null similar curves with variable transformation in $E_1^3$. First, we give the following definition.

**Definition 3.1.** Let $\vec{a}(s_a)$ and $\vec{b}(s_b)$ be Cartan framed null curves in $E_1^3$ with arc length parameters $s_a$, $s_b$, respectively, and let the Frenet frames and invariants of $\vec{a}(s_a)$ and $\vec{b}(s_b)$ be $\{\vec{\alpha}_a, \vec{\gamma}_a, \vec{\beta}_a\}$, $\kappa_a$, $\tau_a$ and $\{\vec{\alpha}_b, \vec{\gamma}_b, \vec{\beta}_b\}$, $\kappa_b$, $\tau_b$, respectively. $\vec{a}(s_a)$ and $\vec{b}(s_b)$ are called null similar curves with variable transformation $\lambda_b^a$ if there exists a variable transformation

$$s_a = \int \lambda_b^a(s_b) ds_b, \tag{10}$$

of the arc lengths such that the null tangents are the same for two curves, i.e.,

$$\vec{\alpha}_a(s_a) = \vec{\alpha}_b(s_b), \tag{11}$$

for all corresponding values of parameters under the transformation $\lambda_b^a$. All null curves satisfying equation (11) are called a family of null similar curves with variable transformation.

If we integrate the (11) we have the following theorem.

**Theorem 3.1.** *The position vectors of the family of Cartan framed null similar curves with variable transformation can be written in the following form*

$$\vec{b}(s_b) = \int \vec{\alpha}_a\left(s_a(s_b)\right) ds_b = \int \vec{\alpha}_a(s_a) \lambda_a^b ds_b.$$

Then we can give the following theorems characterizing null similar curves. In the followings, whenever we talk about Cartan framed null curves $\vec{a}(s_a)$ and $\vec{b}(s_b)$, we mean that the curves have the Frenet frames and invariants as given in Definition 3.1.

**Theorem 3.2.** *Let $\vec{a}(s_a)$ and $\vec{b}(s_b)$ be two null curves in $E_1^3$. Then $\vec{a}(s_a)$ and $\vec{b}(s_b)$ are null similar curves with variable transformation if and only if the principal normal vectors of the curves are the same, i.e.,*

$$\vec{\beta}_a(s_a) = \vec{\beta}_b(s_b), \tag{12}$$

*under the particular variable transformation*

$$\lambda_b^a = \frac{\kappa_b}{\kappa_a}, \tag{13}$$

*of the arc lengths.*

**Proof:** Let $\vec{a}(s_a)$ and $\vec{b}(s_b)$ be two null similar curves in $E_1^3$ with variable transformation. Then differentiating (11) with respect to $s_b$ it follows

$$\kappa_a \vec{\beta}_a \frac{ds_a}{ds_b} = \kappa_b \vec{\beta}_b. \tag{14}$$

From (14), we obtain (12) and (13) immediately.

Conversely, let $\vec{a}(s_a)$ and $\vec{b}(s_b)$ be two null curves in $E_1^3$ satisfying (12) and (13). By multiplying (12) with $\kappa_b$ and differentiating the results equality with respect to $s_b$ we have

$$\int \kappa_b(s_b) \vec{\beta}_b(s_b) ds_b = \int \kappa_b(s_b) \vec{\beta}_b(s_b) \frac{ds_b}{ds_a} ds_a. \tag{15}$$

Then from (12) and (13) we obtain

$$\vec{\alpha}_b(s_b) = \int \kappa_b(s_b) \vec{\beta}_b(s_b) ds_b = \int \kappa_a(s_a) \vec{\beta}_a(s_a) ds_a = \vec{\alpha}_a(s_a), \tag{16}$$

which means that $\vec{a}(s_a)$ and $\vec{b}(s_b)$ are null similar curves with variable transformation.

From Theorem 3.2, we can give the following corollary.

**Corollary 3.1.** *Let $\vec{a}(s_a)$ and $\vec{b}(s_b)$ be two null curves in $E_1^3$. Then $\vec{a}(s_a)$ and $\vec{b}(s_b)$ are null similar curves with variable transformation if and only if $\vec{a}(s_a)$ and $\vec{b}(s_b)$ are null Bertrand curves with particular variable transformation $\lambda_b^a = \frac{\kappa_b}{\kappa_a}$.*

**Theorem 3.3.** Let $\vec{a}(s_a)$ and $\vec{b}(s_b)$ be two null curves in $E_1^3$. Then $\vec{a}(s_a)$ and $\vec{b}(s_b)$ are null similar curves with variable transformation if and only if the binormal vectors of the curves are the same, i.e.,

$$\vec{\gamma}_a(s_a) = \vec{\gamma}_b(s_b), \tag{17}$$

under the particular variable transformation

$$\lambda_b^a = \frac{\tau_b}{\tau_a}, \tag{18}$$

of the arc lengths.

**Proof:** Let $\vec{a}(s_a)$ and $\vec{b}(s_b)$ be two null similar curves in $E_1^3$ with variable transformation. Then, from Definition 3.1 and Theorem 3.2 there exists a variable transformation of the arc lengths such that the tangent vectors and principal normal vectors are the same. Then from (11) and (12) we have

$$\vec{\gamma}_a(s_a) = \vec{\alpha}_a(s_a) \times \vec{\beta}_a(s_a) = \vec{\alpha}_b(s_b) \times \vec{\beta}_b(s_b) = \vec{\gamma}_b(s_b). \tag{19}$$

Conversely, let $\vec{a}(s_a)$ and $\vec{b}(s_b)$ be two null curves in $E_1^3$ satisfying (17) and (18). By differentiating (17) with respect to $s_b$ we obtain

$$\tau_a(s_a)\vec{\beta}_a(s_a)\frac{ds_a}{ds_b} = \tau_b(s_b)\vec{\beta}_b(s_b), \tag{20}$$

which gives us

$$\lambda_b^a = \frac{\tau_b}{\tau_a}, \quad \vec{\beta}_a(s_a) = \vec{\beta}_b(s_b). \tag{21}$$

Then from (17) and (21) we have

$$\vec{\alpha}_a(s_a) = \vec{\beta}_a(s_a) \times \vec{\gamma}_a(s_a) = \vec{\beta}_b(s_b) \times \vec{\gamma}_b(s_b) = \vec{\alpha}_b(s_b) \tag{22}$$

which completes the proof.

**Theorem 3.4.** Let $\vec{a}(s_a)$ and $\vec{b}(s_b)$ be two null curves in $E_1^3$. Then $\vec{a}(s_a)$ and $\vec{b}(s_b)$ are null similar curves with variable transformation if and only if the ratio of curvatures are the same i.e.,

$$\frac{\tau_b(s_b)}{\kappa_b(s_b)} = \frac{\tau_a(s_a)}{\kappa_a(s_a)}, \tag{23}$$

under the particular variable transformation keeping equal total curvatures, i.e.,

$$\varphi_b(s_b) = \int \kappa_b(s_b)ds_b = \int \kappa_a(s_a)ds_a = \varphi_a(s_a) \tag{24}$$

of the arc lengths.

**Proof:** Let $\vec{a}(s_a)$ and $\vec{b}(s_b)$ be two null similar curves in $E_1^3$ with variable transformation. Then from (18) and (21) we have (23) under the variable transformation (24), and this transformation is also leads from (18) by integration.

Conversely, let $\vec{a}(s_a)$ and $\vec{b}(s_b)$ be two null curves in $E_1^3$ satisfying (23) and (24). From Theorem 2.1, the tangents $\vec{\alpha}_a$ and $\vec{\alpha}_b$ of the curves $\vec{a}(s_a)$ and $\vec{b}(s_b)$ satisfy the following vector differential equations of third order

$$\frac{d^3\vec{\alpha}_a}{d\varphi_a^3} + 2f_a(\varphi_a)\frac{d\vec{\alpha}_a}{d\varphi_a} + \frac{df_a(\varphi_a)}{d\varphi_a}\vec{\alpha}_a = 0, \tag{25}$$

$$\frac{d^3\vec{\alpha}_b}{d\varphi_b^3} + 2f_b(\varphi_b)\frac{d\vec{\alpha}_b}{d\varphi_b} + \frac{df_b(\varphi_b)}{d\varphi_b}\vec{\alpha}_b = 0, \tag{26}$$

respectively, where

$$f_a(\varphi_a) = \frac{\tau_a(\varphi_a)}{\kappa_a(\varphi_a)}, \quad f_b(\varphi_b) = \frac{\tau_b(\varphi_b)}{\kappa_b(\varphi_b)}, \quad \varphi_a(s_a) = \int \kappa_a(s_a)ds_a, \quad \varphi_b(s_b) = \int \kappa_b(s_b)ds_b.$$

From (23) we have $f_a(\varphi_a) = f_b(\varphi_b)$ under the variable transformation $\varphi_a = \varphi_b$. So, under the equation (23) and transformation (24), the equations (25) and (26) are the same, i.e., they have the same solutions. It means that the null tangent vectors $\vec{\alpha}_a$ and $\vec{\alpha}_b$ are the same. Then $\vec{a}(s_a)$ and $\vec{b}(s_b)$ are two null similar curves in $E_1^3$ with variable transformation.

Let now consider some special cases. From (13) and (21) we have
$$\kappa_b = \lambda_b^a \kappa_a, \quad \tau_b = \lambda_b^a \tau_a, \tag{27}$$
respectively. From (27) it is clear that if $\vec{a}(s_a)$ is a null geodesic i.e., $\kappa_a = 0$, then under the variable transformation the curvature does not change i.e., $\kappa_b = 0$. Similarly, if $\vec{a}(s_a)$ is a null curve with vanishing torsion i.e., $\tau_a = 0$, then under the variable transformation the torsion $\tau_b$ is given as $\tau_b = 0$. So we have the following corollaries.

***Corollary 3.2.*** *The family of null geodesics forms a family of null similar curves with variable transformation.*

***Corollary 3.3.*** *The family of null curves with vanishing torsion forms a family of null similar curves with variable transformation.*

If $\vec{a}(s_a)$ is a null helix, i.e., $\tau_a/\kappa_a$ is constant, and $\vec{a}(s_a)$ and $\vec{b}(s_b)$ are null similar curves in $E_1^3$ with variable transformation, then from (23) we have $\tau_b/\kappa_b$ is constant. Then we have the following corollary.

***Corollary 3.4.*** *The family of null helices forms a family of null similar curves with variable transformation.*

**4. Conclusions**

A family of null curves in Minkowski 3-space $E_1^3$ are defined and called null similar curves. Some properties of these special curves are obtained and it is showed that null Bertrand curves are null similar curves with variable transformation. Moreover, it is obtained that null geodesics, null curves with vanishing torsion and null helices form the families of null similar curves.